\documentstyle{article}
\title{Coloring finite subsets of uncountable sets}
\author{P. Komj\'{a}th\thanks{Supported by the 
Hungarian OTKA Grant No.T014105}\\E\"{o}tv\"{o}s University 
\and S. Shelah\thanks{No. 516. Research supported by the 
Basic Research Foundation of the 
Israel Academy of Sciences and Humanities}\\Hebrew University}

\newtheorem{thm}{Theorem}
\newtheorem{lem}{Lemma}
\begin{document}
\maketitle
\bibliographystyle{plain}
\def\phi{{\varphi}}
\renewcommand\aa{{\alpha}}
\newcommand\bb{{\beta}}
\renewcommand\gg{{\gamma}}
\newcommand\dd{{\delta}}
\newcommand\ee{{\varepsilon}}
\newcommand\oo{{\omega}}
\def\kk{{\kappa}}
\def\Dom{{\rm Dom}}
\def\Rng{{\rm Ran}}
\def\qed{\hfill{\vbox{\hrule\hbox{\vrule\kern3pt
                \vbox{\kern6pt}\kern3pt\vrule}\hrule}}}
\def\supp{{\rm supp}}
\def\cf{{\rm cf}}
\newcommand\forc{{\,\parallel\joinrel\relbar\joinrel\relbar\,}}

\begin{abstract}
  It is consistent for every \(1\leq n<\oo\) 
  that \(2^\oo=\oo_n\) and there is a function 
   \(F:[\oo_n]^{<\oo}\to\oo\) such that every finite 
   set can be written at most \(2^n-1\) 
   ways as the union of two distinct monocolored sets. 
   If GCH holds, for every such coloring there is a 
   finite set that can be written at least 
   \(\frac{1}{2}\sum^n_{i=1}{n+i\choose n}{n\choose i}\)
   ways as the union of two sets with the same color. 
\end{abstract}
\setcounter{section}{-1}
\section{Introduction}
In \cite{KS} we proved that for every coloring 
\(F:[\oo_n]^{<\oo}\to\oo\) there exists a set 
\(A\in [\oo_n]^{<\oo}\) which can be written at least 
\(2^n-1\) ways as \(A=H_0\cup H_1\) for some 
\(H_0\neq H_1\), \(F(H_0)=F(H_1)\) and that for 
\(n=1\) there is in fact a function \(F\) for which this is 
sharp. 
Here we show that for every \(n<\oo\) it is consistent that 
\(2^\oo=\oo_n\) and for some function $F$ as above 
for every finite set \(A\) there are at most \(2^n-1\) 
solutions of the above equation. 
We use historic forcing which was first used in \cite{BS} and \cite{SS} 
then in \cite{K} and \cite{GS}.
Under GCH, we improve the positive result of \cite{KS} by showing 
that for every  $F$ as above some finite set can be written 
at least 
\(T_n=\frac{1}{2}\sum^n_{i=1}{n+i\choose n}{n\choose i}\) 
ways as the union of two sets with the same $F$ value.

With the methods of \cite{KS} it is easy to show the 
following corollary of our
independence result. 
It is consistent that \(2^\oo=\oo_n\) and 
there is a function \(f:{\bf R}\to\oo\) 
such that if $x$ is a real number then $x$ 
cannot be written more than $2^n-1$ ways as 
the arithmetic mean of some $y\neq z$ with 
\(f(y)=f(z)\). (\((y,z)\) and \((z,y)\) are not regarded 
distinct.) 
Another idea of \cite{KS} can be used to modify 
our second result to the following.
If GCH holds and  $V$ is a vector space over the rationals 
with $|V|=\oo_n$, $f:V\to\oo$ 
then some vector can be written at least $T_n$ ways
as the arithmetic mean of two vectors with the same 
$f$-value. 
\bigskip

\noindent{\bf Notation} 
We use the standard  set theory notation. 
If \(S\) is a set, \(\kappa\) a cardinal, then 
$[S]^\kk=\{A\subseteq S:|A|=\kk\}$, 
$[S]^{<\kk}=\{A\subseteq S:|A|<\kk\}$, 
$[S]^{\leq \kk}=\{A\subseteq S:|A|\leq\kk\}$. 
\(P(S)\) is the power set of \(S\). 
If $f$ is a function, $A$ a set, then 
$f[A]=\{f(x):x\in A\}$.

\section{The independence result}

\begin{thm}
 For\ $1\leq n <\oo$ it is 
 consistent that $2^\oo=\oo_n$ and there is a 
 function $F:[\oo_n]^{<\oo}\to \oo$ such that 
 for every $A\in [\oo_n]^{<\oo}$ there are at most 
 $2^n-1$ solutions of $A=H_0\cup H_1$ with $H_0\neq H_1$, 
 $F(H_0)=F(H_1)$.
\end{thm}
\bigskip

For $\aa<\omega_n$ fix a bijection $\phi_\aa:\aa\to |\aa|$. 
For $x\in [\omega_n]^{<\omega}$ define 
$\gg_i(x)$ for $i<k=\min(n,|x|)$ as follows. 
$\gg_0(x)=\max(x)$. 
$$
  \gg_{i+1}(x)= 
               \phi^{-1}_{\gg_0(x)} \Bigl(\gg_i\bigl(\phi_{\gg_0(x)}
                                        [x\cap \gg_0(x)]\bigr)\Bigr).
$$
$\gg(x)=\{\gg_0(x),\dots,\gg_{k-1}(x)\}$. 

So, for example, if $n=0$ then $\gg(x)=\emptyset$, if $n=1$, 
\(x\neq\emptyset\), then 
$\gg(x)=\{\gg_0(x)\}=\{\max(x)\}$.

\bigskip
\begin{lem}
 Given $s\in [\omega_n]^{\leq n}$ 
 there are at most countably many $x\in [\omega_n]^{<\omega}$ 
 such that $\gg(x)=s$.
\end{lem}
\bigskip
\noindent{\bf Proof} 
By induction on $n$. 
\qed
\bigskip

Let $\Phi(s)=\bigcup\{x:\gg(x)\subseteq s\}$, 
a countable set for \(s\in [\oo_n]^{<\oo}\).

\bigskip
\noindent{\bf Definition} The two sets 
$x$, $y\in [\omega_n]^{<\omega}$ are 
{\sl isomorphic} if the structures 
$(x;<,\gg_0(x), \dots, \gg_{k-1}(x))$, 
$(y;<,\gg_0(y), \dots, \gg_{k-1}(y))$, are isomorphic, i.e., 
$|x|=|y|$ and the positions of the elements $\gg_i(x)$, $\gg_i(y)$ 
are the same. 

\medskip
Notice that for every finite $j$ there are just finitely many 
isomorphism types of $j$-element sets. 
\medskip

The elements of $P$, the applied notion of forcing will be 
{\sl some} structures of the form 
$p=(s,f)$ where $s\in[\omega_n]^{<\omega}$ and 
$f:P(s)\to \omega$.

The only element of $P_0$ is ${\bf 1}_P=(\emptyset,\langle \emptyset,0
\rangle)$,  it will be the largest element of $P$. 
The elements of $P_1$ are of the form 
$p=(\{\xi\},f)$ where $f(\emptyset)=0\neq f(\{\xi\})$ 
for $\xi<\omega_n$. 

Given $P_t$, 
$p=(s,f)$ is in  $P_{t+1}$ if the following is true. 
$s=\Delta\cup a\cup b$ is a disjoint decomposition. 
$p'=(\Delta\cup a,f')$ and $p''=(\Delta\cup b,f'')$ are in $P_t$ 
where $f'=f|P(\Delta\cup a)$, $f''=f|P(\Delta\cup b)$. 
There  is $\pi:\Delta\cup a\to \Delta\cup b$, an isomorphism 
between $(\Delta\cup a,<,P(\Delta\cup a),f')$ and 
$(\Delta\cup b,<,P(\Delta\cup b),f')$. 
$\pi|\Delta$ is the identity. 
For $H\subseteq \Delta\cup a$ the sets $H$ and 
$\pi[H]$ are isomorphic. 
$a\cap \Phi(\Delta)=b\cap \Phi(\Delta)=\emptyset$. 
$f-f'-f''$ is one-to-one and takes only values outside 
$\Rng(f')$ (which is the same as $\Rng(f'')$). 
$P=\bigcup\{P_t:t<\omega\}$. 
We make $p\leq p',p''$ and the ordering on $P$ is the one 
generated by this. 
\bigskip
\begin{lem}
 $(P,\leq)$ is ccc.
\end{lem}
\bigskip
\noindent{\bf Proof} 
Assume that $p_\aa\in P$ ($\aa<\omega_1$). 
We can assume by thinning and using the $\Delta$-system lemma 
and the pigeon hole principle that the following hold. 
$p_\aa\in P_t$ for the same $t<\omega$. 
$p_\aa=(\Delta\cup a_\aa,<,P(\Delta\cup a_\aa),f_\aa)$ where the structures 
$(\Delta\cup a_\aa,<,f_\aa)$ and 
$(\Delta\cup a_\bb,<,f_\bb)$ are isomorphic for $\aa,\bb<\oo_1$, 
$\{\Delta,a_\aa:\aa<\oo_1\}$ pairwise disjoint. 
We can also assume that if $\pi$ is the 
isomorphism between 
$(\Delta\cup a_\aa,<,f_\aa)$ and 
$(\Delta\cup a_\bb,<,f_\bb)$ then $H$ and $\pi[H]$ 
are isomorphic for $H\subseteq \Delta\cup a_\aa$. 
Moreover, if we assume that $\Delta$ occupies the 
same positions in the ordered sets $\Delta\cup a_\aa$ 
($\aa<\oo_1$) then $\pi$ will be the identity on $\Delta$. 
As $\Phi(\Delta)$ is countable, by removing 
countably many indices we can also assume that 
$\Phi(\Delta)\cap a_\aa=\emptyset$ for $\aa<\oo_1$. 
Now any $p_\aa$ and $p_\bb$ are compatible as we can 
take 
$p=(\Delta\cup a_\aa\cup a_\bb,<,P(\Delta\cup a_\aa \cup a_\bb),f)
\leq p_\aa$, $p_\bb$ 
where $f\supseteq f_\aa$, $f_\bb$ 
is an appropriate extension, i.e., 
$f-f_\aa -f_\bb$ is one-to-one and takes values 
outside $\Rng(f_\aa)$. 
\qed

\bigskip
\begin{lem}
 If $(s,f)\in P$, $H_0,H_1\subseteq s$ 
 have $f(H_0)=f(H_1)$ then $H_0$, $H_1$ are isomorphic.
\end{lem}
\bigskip
\noindent{\bf Proof} 
Set  $(s,f)\in P_t$. 
We prove the statement by induction on $t$. 
There is nothing to prove for $t<2$. 
Assume now that $(s,f)\in P_{t+1}$, 
$s=\Delta\cup a\cup b$, 
$\pi:\Delta\cup a\to \Delta\cup b$ 
as in the definition of $(P,\leq)$. 
As $f(H_0)$ is a value taken twice by $f$, both 
$H_0$ and $H_1$ must be subsets of either 
$\Delta\cup a$ or $\Delta\cup b$. 
We are done by induction unless $H_0\subseteq \Delta\cup a$ 
and $H_1\subseteq \Delta\cup b$ (or vice versa). 
Now $H_0$ and $\pi[H_0]$ are isomorphic and 
$f(H_0)=f(\pi[H_0])=f(H_1)$ so by the inductive hypothesis 
$\pi[H_0]$ and $H_1$ are ismorphic and then so are $H_0$, $H_1$. 
\qed

\bigskip
\begin{lem}
 If $(s,f)\in P$, 
 $H_0,H_1\subseteq s$, 
 $f(H_0)=f(H_1)$, $x\in H_0\cap H_1$ then $x$ occupies the 
 same position in the ordered sets $H_0$, $H_1$.
\end{lem}
\bigskip
\noindent{\bf Proof} 
Similarly to the proof of the previous Lemma, 
by induction on $t$, for $(s,f)\in P_t$. 
With similar steps, we can assume that 
$(s,f)=(\Delta\cup a\cup b,f)\leq (\Delta\cup a,f'), (\Delta\cup b,f'')$, 
$H_0\subseteq \Delta\cup a$, 
$H_1\subseteq \Delta\cup b$. 
Notice that $x\in \Delta$. 
Now, as $\pi(x)=x$, 
$x$ is a common element of $\pi[H_0]$ and $H_1$ and also 
$f''(\pi[H_0])=f''(H_1)$. 
By induction we get that $x$ occupies the same position in 
$\pi[H_0]$ and $H_1$ so by pulling back we get that this is 
true for $H_0$ and $H_1$. 
\qed

\bigskip
\begin{lem}
 If $(s,f)\in P$, 
 $A\subseteq s$, $0\leq j\leq n$ then 
 $A$ can be written at most $2^j-1$ ways as 
 $A=H_0\cup H_1$ with $H_0$, $H_1$ distinct, 
 $f(H_0)=f(H_1)$, and $|\gg(H_0)\cap \gg(H_1)|\geq n-j$.
\end{lem} 
\bigskip
\noindent{\bf Proof} 
By induction on $j$ and inside that induction, by induction on $t$, 
for $(s,f)\in P_t$. 
The case $t<2$ will always be trivial. 

Assume first that $j=0$. 
In this case our Lemma reduces to the following statement. 
There are no $H_0\neq H_1$ such that $\gg(H_0)= \gg(H_1)$. 
In the inductive argument we assume as usual that 
$s=\Delta\cup a\cup b$ and so $(s,f)\in P_{t+1}$ was created from 
$(\Delta\cup a,f')$ and $(\Delta\cup b,f'')$, 
$H_0\subseteq \Delta\cup a$, 
$H_1\subseteq \Delta\cup b$. 
As 
$\gg(H_0)= \gg(H_1)$, 
$\gg(H_0)\subseteq \Delta$, but then, 
as $\Phi(\Delta)\cap a=\emptyset$, $H_0$ can have no points outside 
$\Delta$ and similarly for $H_1$, so we can go back, say to 
$(\Delta\cup a,f')\in P_t$ which concludes the argument. 

Assume now that the statement is proved for 
$j$ and we have $p=(s,f)\in P_{t+1}$, 
$s=\Delta\cup a\cup b$ and  $p$ was created from 
$p'=(\Delta\cup a,f')$ and $p''=(\Delta\cup b,f'')$. 
In $A\subseteq \Delta\cup a\cup b$ we can assume that 
$y=A\cap a \neq\emptyset$, $z=A\cap b\neq \emptyset$ as 
otherwise we can pull back to $p'$ or $p''$. 
But then, if $A=H_0\cup H_1$, then,if,  say, 
$H_0 \subseteq \Delta\cup a$, $H_1\subseteq \Delta\cup b$ hold, 
then necessarily $H_0\cap a=y$, $H_1\cap b=z$, 
so $H_0=x_0\cup y$, $H_1=x_1\cup z$ where 
$x_0\cup x_1=x=A\cap \Delta$. 
We can create decompositions of 
$B=x\cup\pi[y]\cup z$ by taking $B=\pi[H_0]\cup H_1$. 
But some of these decompositions will not be different and 
it may happen that we get non-proper (i.e., one-piece) decomposition. 
This can only happen if $\pi[y]=z$, and then 
the two decompositions 
$A=(x_0\cup y)\cup(x_1\cup z)$ and $A=(x_1\cup y)\cup(x_0\cup y)$ 
produce the same decomposition of $B$, namely, 
$B=(x_0\cup z)\cup(x_1\cup z)$ and there is but one decomposition, 
$A=(x\cup y)\cup(x\cup z)$ which cannot be mapped to a decomposition 
of $B$. If this (i.e., \(\pi[y]=z\)) does not happen, 
we are done by induction. 
If this does happen, we know that $\gg(H_0)=\gg(x_0\cup y)$ has 
an element in $y$ (by the argument at the beginning of the proof). 
As $f(x_0\cup y)=f(x_1\cup z)$, by Lemmas 3 and 4, both 
$H_0=x_0\cup y$ and $H_1=x_1\cup z$ have an element in the $\gg$-subset, 
at the same positions  which are mapped onto each other by $\pi$. 
We get that $\gg(x_0\cup z)\cap\gg(x_1\cup z)$ has at least $n-j$ 
element, so by our inductive assumption we have at most 
$2^j-1$ decompositions, which gives at most 
$2\cdot(2^j-1)+1=2^{j+1}-1$ decompositions of $A$. 

\qed

\bigskip
Let $G\subseteq P$ be a generic subset. 
Set $S=\bigcup\{s:(s,f)\in G\}$,  
$F=\bigcup\{f:(s,f)\in G\}$. 

\bigskip
\begin{lem}
 There is a $p\in P$ such that 
 $p\forc |S|=\aleph_n$.
\end{lem}
\bigskip
\noindent{\bf Proof} 
Otherwise ${\bf 1}\forc \sup(S)<\oo_n$. 
By ccc, there is an ordinal $\xi<\oo_n$ for which  
${\bf 1}\forc \sup(S)<\xi$, but this is impossible  as there are 
conditions in $P_1$ forcing that $\xi\in S$.  
\qed

\bigskip
Now we can conclude the proof of the Theorem. 
If $G$ is generic, and $p\in G$ with the condition $p$ of Lemma 6, 
then in $V[G]$   $F$ witnesses the theorem 
by Lemma 5 (for \(j=n\)) on the ground set 
\(S\). 
As \(|S|=\oo_n\) we can replace it by \(\oo_n\). \qed

\section{The GCH result}

Set 
$$
T_n= \frac{1}{2}\sum^n_{i=1}{n+i\choose n}{n\choose i}.
$$
So \(T_1=1\), \(T_2=6\), 
\(T_3=31\). 
In general, \(T_n\) is asymptotically 
\(c(3+2\sqrt{2})^n/\sqrt{n}\)  for some $c$.
\begin{thm}
 {\em (GCH)} If $F:[\oo_n]^{<\oo}\to \oo$ 
 then some $A\in [\oo_n]^{<\oo}$ has 
 at least $T_n$ decompositions as 
 $A=H_0\cup H_1$, $H_0\neq H_1$, 
 $F(H_0)=F(H_1)$.
\end{thm}
\bigskip
\noindent{\bf Proof} 
By the Erd\H{o}s-Rado theorem (see \cite{EHMR,ER}) 
there is a set $\{x_\aa:\aa<\oo_1\}$ which is 
$(n-1)$-end-homogeneous, i.e., for some 
$g:[\oo_1]^{<\oo}\to\oo$, 
if 
$\aa_1<\cdots<\aa_k<\bb_1<\cdots<\bb_{n-1}<\oo_1$ 
then 

$$
f(\{ x_{\aa_1},\dots,x_{\aa_{k}},x_{\bb_1},\dots,x_{\bb_{n-1}}\}) 
= g(\aa_1,\dots,\aa_k).
$$

Select 
$S_1\in [\oo_1]^{\oo_1}$ 
in such a way that $g(\aa)=c_0$ for $\aa\in S_1$. 
Set $\gg_1=\min(S_1)$. 
In general, if $\gg_i$, $S_i$ are given ($1\leq i<n$) 
pick 
$S_{i+1}\in [S_i-(\gg_i+1)]^{\oo_1}$ 
so that 
$g(\gg_1,\dots,\gg_i,\aa)=c_i$ 
for $\aa\in S_{i+1}$ and set 
$\gg_{i+1}=\min(S_{i+1})$. 
Given $\gg_1,\dots,\gg_n$ and 
$S_n$ let 
$\gg_{n+1},\dots,\gg_{2n}$ be the $n$ 
least elements 
of $S_n-(\gg_n+1)$. 

Our set will be 
$A=\{x_{\gg_1},\dots,x_{\gg_{2n}}\}$. 
For $0\leq i<n$ the color of any $(n+i)$-element 
subset of $A$ containing $x_{\gg_1},\dots,x_{\gg_i}$ will be $c_i$. 
We can select 
$\frac{1}{2}{2n-i\choose n}{n\choose i}$ 
different pairs of those sets which cover $A$. 
In toto, we get $T_n$ decompositions of $A$. \qed
%


\bigskip

\hbox{\vbox{\hbox{P\'eter Komj\'ath}
            \hbox{Department of Computer Science}
            \hbox{E\"otv\"os University}
            \hbox{Budapest, M\'uzeum krt.~6--8}
            \hbox{1088, Hungary}
            \hbox{e-mail: {\tt kope@cs.elte.hu}}
     }
     \hskip3cm
     \vbox{\hbox{Saharon Shelah}
           \hbox{Institute of Mathematics}
           \hbox{Hebrew University}
           \hbox{Givat Ram, 91904}
           \hbox{Jerusalem, Israel}
           \hbox{e-mail: {\tt shelah@math.huji.ac.il}}}
      }


\begin{thebibliography}{EHMR}
\bibitem{BS} J. E. Baumgartner, S. Shelah. 
   Remarks on superatomic Boolean algebras, 
   {\sl Annal of Pure and Applied Logic \bf 33} (1987),
   109--129.

\bibitem{EHMR} P. Erd\H os, A. Hajnal, A. M\'at\'e, R. Rado. 
   {\sl Combinatorial Set Theory: Partition Relations for Cardinals}, 
   North-Holland, Studies in Logic, {\bf 106}, (1984).

\bibitem{ER} P.~Erd\H os, R. Rado. 
   A partition calculus in set theory, 
   {\sl Bull. Amer. Math. Soc. \bf 62} (1956), 427--289.

\bibitem{GS} M.~Gilchrist, S.~Shelah. 
   On identities of colorings of pairs for $\aleph_n$, submitted.

\bibitem{K} P.~Komj\'ath. 
   A set mapping with no infinite free subsets, 
   {\sl Journal of Symbolic Logic \bf 56} (1991), 1400--1402.

\bibitem{KS} P.~Komj\'ath, S.~Shelah.
   On uniformly antisymmetric functions, 
   {\sl Real Analysis Exchange \bf   19} (1993--1994), 218--225.

\bibitem{SS} S. Shelah, L. Stanley:
   A theorem and some consistency results in partition calculus, 
   {\sl Annal of Pure and Applied Logic \bf 36} (1987),
   119--152.
\end{thebibliography}
\end{document}